\input{aipcheck}

\documentclass[
    ,final            
  ]
  {aipproc}

\layoutstyle{6x9}

\usepackage{amssymb}
\usepackage{amsfonts}
\usepackage{latexsym}
\usepackage{amsxtra, amsmath}
\usepackage{verbatim}

\newtheorem{theorem}{Theorem}

\newtheorem{proposition}[theorem]{Proposition}

\newtheorem{definition}{Definition}
\newtheorem{example}{Example}

\newtheorem{remark}{Remark}

 \newcommand{\nc}{\newcommand}

\renewcommand{\aa}{\mathfrak{a} }
\nc{\dd}{\mathfrak{d} } 
\nc{\ggo}{\mathfrak{g} }
\nc{\hh}{\mathfrak{h} }
\nc{\mm}{\mathfrak{m} }   
\nc{\nn}{\mathfrak{n} }
\nc{\sso}{\mathfrak{so} }  
\nc{\vv}{\mathfrak{v} } 
\nc{\zz}{\mathfrak{z} }

\nc{\NN}{{\mathbb N}}
\nc{\RR}{{\mathbb R}}  
\nc{\ZZ}{{\mathbb Z}}

\nc{\glg}{\mathfrak{gl} }
\nc{\la}{\langle} \nc{\ra}{\rangle}
\nc{\brg}{[\,,\,]_{\ggo}}
\nc{\brv}{[\,,\,]_{\vv}}
 
 \nc{\SO}{{\sf SO}} \nc{\Spe}{{\sf Sp}} \nc{\Sl}{{\sf Sl}}
 \nc{\SU}{{\sf SU}} \nc{\Or}{{\sf O}} \nc{\U}{{\sf U}}
 \nc{\Gl}{{\sf Gl}} \nc{\Se}{{\sf S}} \nc{\Cl}{{\sf Cl}}
 \nc{\Spin}{{\sf Spin}} \nc{\Pin}{{\sf Pin}}

 \nc{\ad}{\operatorname{ad}} 
 \nc{\Ad}{\operatorname{Ad}}
 
\nc{\Dera}{\operatorname{Dera}}
\nc{\tr}{\operatorname{tr}}
\nc{\ct}{\operatorname{T}}

\nc{\mcn}{\multicolumn}
\nc{\ovr}{\overset{}}
\nc{\undr}{\underset{}}
\nc{\f}{\frac}

\begin{document}

\title{Examples of naturally reductive pseudo Riemannian Lie groups}



\author{Gabriela P. Ovando}{
  address={CONICET - Departamento de Matematica, ECEN - FCEIA, Universidad Nacional de Rosario, Pellegrini 250, 2000 Rosario, Argentina, Email: gabriela@fceia.unr.edu.ar}
}


\begin{abstract}
We provide examples  of naturally reductive pseudo-Riemannian spaces, in particular an example of a  naturally reductive pseudo-Riemannian  2-step 
nilpotent Lie group $(N, \la \,,\,\ra_N)$, such that $\la \,,\,\ra_N$ is
invariant under a left action and for which the center is degenerate. The metric does not correspond to a bi-invariant one.
\end{abstract}

\maketitle


\section{Introduction}

Recent advances in mathematics and physics have renovated the interest in g.o. spaces, that is, manifolds in which every geodesic is homogeneous (see \cite{Du}). 
A subclass of g.o. spaces is provided by the naturally reductive pseudo-Riemannian manifolds.

\begin{definition} A homogeneous pseudo-Riemannian manifold $M$
is said to be {\em naturally reductive} if there is a transitive Lie group of isometries
$G$ with Lie algebra $\ggo$  and there exists a subspace
$\mm\subseteq \ggo$ complementary to $\hh$ in $\ggo$,  $\hh$ the Lie algebra of the isotropy group 
$H$,  such that
$$\Ad(H)\mm \subseteq \mm \quad \mbox{ and }\quad 
\la [x,y]_{\mm}, z\ra + \la y, [x,z]_{\mm} \ra=0 \qquad \mbox{ for all } x, y,
 z\in \mm.$$
\end{definition}
Frequently we will say that a metric on a homogeneous space $M$ is naturally reductive even though it is not naturally reductive with respect to a particular transitive group of isometries (see Lemma 2.3 in \cite{Go}).

  Indeed pseudo-Riemannian
 symmetric spaces are naturally reductive.  Other examples of naturally
 reductive spaces arise from Lie groups equipped with a bi-invariant metric,
 which could exist for nilpotent ones (see \cite{F-S, M-R}). 
 
 In the Riemannian case Gordon proved that a naturally reductive Riemannian nilmanifold can be at most 2-step nilpotent. In \cite{Ov2} the author investigated  2-step nilpotent Lie groups equipped with a naturally reductive pseudo-Riemannian metric. In particular, for the case of nondegenerate center,  necessary and sufficient conditions are given for a indefinite metric to be naturally reductive. In the case of degenerate center, the known examples are provided by the bi-invariant metrics.
 
 In \cite{Ko} Kostant proved that if $M$ denotes a naturally reductive Riemannian space, such that the action of the isometry group $G$ on $M$ is transitive and almost effective, then $G$ can be provided with a bi-invariant metric. In \cite{Ov3} this  result was revalidated for the pseudo-Riemannian case. As a consequence, naturally reductive
metrics can be produced in the following way. Let $\ggo$ denote a Lie algebra equipped with an ad-invariant metric $Q$ and let
$\hh \subset \ggo$ be a nondegenerate Lie subalgebra. Thus one has the following reductive decomposition
\begin{equation}\label{e1}
\ggo = \hh \oplus \mm\qquad \mbox{ with } [h,m] \subseteq \mm \quad \mbox{ and } \quad  m = h^{\perp}.
\end{equation}
Let $G$ denote a Lie group with Lie algebra $\ggo$ and endowed with the bi-invariant metric induced by $Q$ and let $H \subset G$ denote a
closed Lie subgroup with Lie algebra $\hh$. Then the coset space $G/H$ becomes a naturally reductive pseudo-Riemannian space. 

The splitting given in (\ref{e1}) can be modified in some cases in order to produce a naturally reductive pseudo-Riemannian Lie group (see \cite{Ov3} for more details). Here we shall make use of this method in order to see an example of a naturally reductive (not bi-invariant) metric on a 2-step nilpotent Lie group for which  the center is degenerate.

 \section{Naturally reductive compact examples}

 An {\em ad-invariant metric} on a Lie algebra $\ggo$ is a nondegenerate symmetric
 bilinear map $\la\,,\,\ra:\ggo \times \ggo \to \RR$ such that 
\begin{equation}
\la [x,y], z\ra + \la y, [x,z]\ra =0 \qquad \mbox{ for all } x,y, z \in \nn.
\end{equation}
Recall that on a connected Lie group $G$ furnished with a left-invariant
pseudo-Riemannian metric $\la\,,\,\ra$ the following  statements are equivalent
(see \cite{ON} Ch. 11):

\begin{enumerate}
\item $\la\,,\,\ra$ is right invariant, hence bi-invariant;
\item $\la\,,\,\ra$ is $\Ad(G)$-invariant;
\item the inversion map $g\to g^{-1}$ is an isometry of $G$;
\item $\la [x,y], z\ra + \la y, [x,z]\ra =0$ for all $x,y, z \in \ggo$;
\item $\nabla_xy =\frac12 [x,y]$ for all $x,y\in \ggo$, where $\nabla$ denotes the
Levi Civita connection;
\item the geodesics of $G$ starting at $e$ are the one parameter subgroups of $G$.
\end{enumerate}

 Clearly $(G, \la\,,\,\ra)$ is naturally
reductive, which by (3)   is a symmetric space.  Furthermore by computing the curvature tensor one has 
$$R(x,y)=-\frac14 \ad([x,y])\qquad \mbox{ for }x,y \in \ggo.$$

Hence any  {\em simply connected 2-step nilpotent Lie group equipped with a 
bi-invariant  metric is flat}.

The set of nilpotent Lie groups carrying a bi-invariant pseudo-Riemannian metric is  non empty (see \cite{F-S}).  Otherwise  in the Riemannian case, a naturally reductive  nilpotent Lie group may be at most 2-step nilpotent \cite{Go}.

Among other possible constructions, 2-step nilpotent Lie algebras admitting an
ad-invariant metric can be obtained as follows. Let $(\vv, \la \,,\, \ra_+)$ denote
a real vector space equipped with an inner product and let $\rho:\vv \to \sso(\vv,\la\,,\,\ra_+)$
an injective linear map satisfying 
\begin{equation}\label{jad}
\rho(u)u=0\qquad\mbox{  for all }u\in \vv.
\end{equation}
 Consider the
vector space $\nn:=\vv^*\oplus \vv$ furnished with the canonical neutral metric
$\la\,,\,\ra$ and
define a Lie bracket on $\nn$ by
\begin{equation}\label{brad}
\begin{array}{rcl}
[x,y] & = & 0 \quad \mbox{ for  }x\in \vv^*, y\in \nn\quad \mbox{ and }\quad 
[\nn,\nn]\subseteq \vv^*\\
\la [u,v], w\ra & = & \la \rho(w) u,v\ra_+ \qquad \mbox{ for all } u,v,w\in \vv.
\end{array}
\end{equation}
 Then $\nn$ becomes a 2-step nilpotent Lie algebra of
corank zero for which the metric $\la\,,\,\ra$ is ad-invariant. This construction was called the
{\em modified cotangent}, since $\nn$ is linear isomorphic to the cotangent of $\vv$.
Notice that the commutator coincides with the center and it equals $\vv^*$. This
allows to construct 2-step nilpotent Lie algebras of null corank which carry an ad-invariant metric.  Furthermore this is basically
the way to obtain  such Lie algebras, (see \cite{N-R, Ov1}  for more
details):

\begin{theorem} \label{mod}  Let $(\nn, \la\,,\,\ra)$ denote a 2-step nilpotent Lie algebra of corank
$m$ endowed with an ad-invariant metric. Then $(\nn, \la\,,\,\ra)$ is isometric
isomorphic to an orthogonal direct product of the Lie algebras $\RR^m$ and a
modified cotangent.
\end{theorem}

 Examples of 2-step nilpotent Lie algebras with ad-invariant metrics arise
  by taking $\ct^*\nn$, the cotangent of any 2-step
nilpotent Lie algebra $\nn$ together with the canonical neutral metric. Let $\nn=\zz\oplus \vv$ denote a 2-step nilpotent Lie
algebra, where $\vv$ is any complementary subspace of $\zz$ in $\ggo$. Let $z_1, \hdots, z_m$ be a basis of the center $\zz$ and let $v_1,
\hdots, v_n$  be a  basis of the vector space $\vv$. Thus
$$[v_i, v_j]=\sum_{s=k}^m c_{ij}^s z_s \qquad \quad i, j=1, \hdots n.$$ 
Let $\ct^*\nn=\nn\ltimes \nn^*$ denote the cotangent Lie algebra obtained via
the coadjoint representation. Indeed  the set
 $\{z^1, \hdots, z^m, v^1,
\hdots, v^n\}$ becomes the  dual basis of the basis above adapted to the
decomposition $\nn^*=\zz^*\oplus \vv^*$. The non trivial Lie bracket relations concerning
the coadjoint action follow
$$
[v_i, z^j]=\sum_{s=1}^n d_{ij}^s v^s \qquad \quad \mbox{ for } i=1, \hdots n, j=1, \hdots m.
$$
Thus 
$[v_i, z^j](v_k)= d_{ij}^k$ and by the definition 
$$
[v_i, z^j](v_k) =-z^j(\sum_{s=1}^m c_{ik}^s z^s)= - c_{ik}^j \qquad \quad 
i,k=1, \hdots n, j=1, \hdots m.
$$
Therefore $d_{ij}^k= - c_{ik}^j$ for $i,k=1, \hdots n, j=1, \hdots m$.

It is clear that if for some basis of $\nn$ the structure constants are rational numbers then by choosing the union of this basis and its dual on $\ct^*\nn$ one gets
rational structure constants for $\ct^*\nn$. Thus by the Mal'cev criterium $N$
and  its cotangent $\ct^*N$, the simply connected Lie group with Lie algebra $\ct^*\nn$, admits
a lattice  which induces   a compact quotient  (see \cite{O-V,Ra} for instance).
 
Let $\Gamma\subset \ct^*N$ denote a cocompact subgroup of $\ct^*N$. Indeed
$\ct^*N$ acts on the compact nilmanifold $(\ct^*N)/\Gamma$ by left translation isometries if we
induce to the quotient the bi-invariant metric corresponding to the neutral canonical one on
$\ct^*\nn$. The tangent
space at the representative $e$ can be identified with $\ct^*\nn \simeq T_e((\ct^*N)/\Gamma)$ so that
$\ct^*\nn=\{0\}\oplus \ct^*\nn$ and clearly $Ad(\Gamma)\ct^*\nn\subseteq
\ct^*\nn$ which says that $(\ct^*N)/\Gamma$ is homogeneous reductive (see
\textsection 3 Ch. X vol. 2 \cite{K-N}). Moreover 
 the induced  metric on the quotient  satisfies
$$\la[x,y],z\ra+ \la [x,z], y\ra=0\qquad \quad \forall x,y, z\in \ct^*\nn.$$

\begin{proposition} Let $N$ denote a 2-step nilpotent Lie group. If it admits a
rational lattice then  the cotangent Lie group $\ct^*N$ admits a cocompact
subgroup $\Gamma$ such that $(\ct^*N)/\Gamma$ is a naturally reductive pseudo-Riemannian space.
\end{proposition}

\begin{example} The low dimensional 2-step nilpotent Lie group  admitting an 
ad-invariant metric occurs in dimension six. 
 This Lie algebra is  the cotangent of the Heisenberg Lie algebra $\ct^*\hh_3$. Explicitly let $e_1,e_2,e_3, e_4, e_5, e_6$ be a basis of $\ct^*\hh_3$; the Lie brackets are
$$[e_4,e_5]=e_1 \qquad [e_4,e_6]=e_2 \qquad [e_5,e_6]=e_3$$
and an ad-invariant metric $Q$ is defined by the non zero symmetric relations
\begin{equation}\label{e5}
1 = \la e_1, e_6\ra =\la e_2, e_5\ra =\la e_3, e_4\ra.
\end{equation}

Let $\ct^* H_3$ denote the corresponding simply connected six dimensional Lie group with Lie algebra $\ct^*\hh_3$. This  can be modelled on $\RR^6$ together with  the multiplication group  given by
$$\begin{array}{rcl}
(x_1,x_2,x_3, x_4,x_5,x_6) \cdot (y_1, y_2,y_3,y_4,y_5,y_6) &  = & (x_1 + y_1 + \frac12(x_4y_5-x_5y_4), \\
&& x_2+y_2+\frac12(x_4y_6-x_6y_4),\\
& &  x_3+y_3+\frac12(x_5y_6-x_6y_5),\\
&& x_4+y_4, x_5+y_5, x_6+y_6).
\end{array}
$$
By the Malcev criterium $\ct^* H_3$ admits a cocompact lattice $\Gamma$. 
 By inducing the bi-invariant metric of $\ct^* H_3$ to $\ct^* H_3/\Gamma$ one gets a 
 invariant metric on  $\ct^* H_3/\Gamma$, and in this way $\ct^* H_3/\Gamma$ is
 a pseudo-Riemannian naturally reductive compact nilmanifold.

For instance the subgroup of $\ct^* H_3$ given by
$$\Gamma=\{(k_1,k_2,k_3, 2k_4, k_5, 2k_6)\,/\mbox{ for }\,k_i\in \ZZ \, \forall i=1,2,3,4,5,6\}$$
is a co-compact lattice of $\ct^* H_3$, so that $\ct^* H_3/\Gamma$ is a 
compact homogeneous manifold.
\end{example}

\section{Naturally reductive noncompact examples}

The family of naturally reductive pseudo-Riemannian Lie groups constructed in \cite{Ov3} is obtained as follows. Take
\begin{itemize}
\item a Lie algebra $\hh$ with ad-invariant metric $\la \cdot , \cdot \ra_{\hh}$,

\item a Lie algebra $\dd$ with ad-invariant metric $\la \cdot , \cdot \ra_{\dd}$,

\item a Lie algebra homomorphism $\pi: \hh \to \Dera(\dd, \la \cdot , \cdot \ra_{\dd})$ from $\hh$ to the Lie algebra of skew-symmetric derivations of
$(\dd, \la \cdot , \cdot \ra_{\dd})$.
\end{itemize}

This data gives rise to a Lie algebra $\ggo$ which can be provided with an ad-invariant metric $Q$. Moreover  $\ggo$ admits  a splitting as in (\ref{e1}), so that $\hh$ and $\mm=\hh^{\perp}$ are nondegenerate with respect to $Q$ and $\ggo$ also decomposes  as a semidirect sum of vector spaces as follows
\begin{equation}\label{e2}
\ggo = \hh \oplus \mathcal G(\dd),
\end{equation}
where $\mathcal G(\dd)$ denotes the Lie algebra of the Lie group $\mathcal G(D)$ and $\ggo$ is an isometry Lie algebra acting on $\mathcal G(\dd)$ with the stability Lie algebra $\hh$. In
general $\mm$ as in (\ref{e1}) does not coincide with $\mathcal G(\dd)$, however as vector spaces they are isomorphic via the map $\lambda: \mathcal G(\dd) \to \mm$ which induces the metric of $\mm$ to $\mathcal G(\dd)$ making both spaces linearly isometric. The metric $\la \cdot , \cdot \ra$ induced on $\mathcal G(\dd)$ is defined on the Lie group $\mathcal G(D)$ by translations on the left.

Below we shall  make use of this construction to provide an example of a naturally reductive pseudo Riemannian 2-step nilpotent Lie group, such that its center is degenerate.

\medspace

Let  $\aa(1,2)$ denote the free 3-step nilpotent Lie algebra with basis $\{e_1,e_2,e_3,e_4,e_5\}$ and the  non-zero Lie brackets 
$$[e_1,e_2]=e_3,\qquad [e_1,e_3]=e_4,\qquad [e_2,e_3]=e_5. $$
Endow $\aa(1,2)$ with the ad-invariant metric  $Q$  given by
\begin{equation}\label{e4}
-Q(e_1, e_1)= Q(e_1, e_5)= - Q(e_2, e_4)= Q(e_3, e_3)=1.
\end{equation}

Indeed with respect $Q$, $\hh= \RR e_1$ is a nondegenerate subalgebra of $\aa(1,2)$ such that $\hh^{\perp}$ is 
$$\hh^{\perp}=\RR (e_1+e_5) \oplus \RR e_2 \oplus \RR e_3 \oplus \RR e_4.$$

Let $G$ denote the nilpotent Lie group with Lie algebra $\aa(1,2)$ and let $H\subset G$ the closed Lie subgroup with Lie algebra $\hh$. According to \cite{Ko,Ov3} the homogeneous space $G/H$ becomes a naturally reductive space. Now we follow the construction of \cite{Ov3} to produce a naturally reductive metric on the Lie group $\mathcal G(D)$ whose Lie algebra is $\RR \oplus \hh_3 = span\{e_2, e_3, e_4, e_5\}$. Notice that $\nn:=span\{e_2, e_3, e_4, e_5\}$ is a degenerate ideal in $\aa(1,2)$, which is a 2-step nilpotent Lie algebra of dimension four isomorphic to the trivial extension of the Heisenberg Lie algebra of dimension three $\hh_3=span\{e_2, e_3, e_5\}$.

Thus $e_1$ acts on $\RR \oplus \hh_3$ as a nilpotent derivation
$$\left( 
\begin{matrix}
0 & 0 & 0 & 0\\
1 & 0 & 0 & 0\\
0 & 1 & 0 & 0 \\
0 & 0 & 0 & 0
\end{matrix}
\right)
$$

Induce on $\nn$ the following metric
\begin{equation}\label{e3}
- \la e_2, e_4\ra = \la e_3, e_3 \ra = \la e_5, e_5\ra=1.
\end{equation}

With respect to this metric the center of $\nn$ which is spanned by $e_4, e_5$ is degenerate and Theorem 3.1 in \cite{Ov3} shows that this metric gives rise to a naturally reductive metric on the corresponding Lie group of $\nn$. This metric is not ad-invariant. See for instance \cite{Ov1}.

\begin{remark} The Lie algebra $\RR \oplus \hh_3$ is related to the space known as the Kodaira Thurston manifold, which provides an example of a manifold admitting symplectic but not K\"ahler structures.
\end{remark}

\subsection{Naturally reductive pseudo-Riemannian  nilmanifolds in low dimensions}

\begin{enumerate}
\item dimension three: $H_3$, the Heisenberg Lie group equipped with pseudo Riemannian metric, such that the center is nondegenerate \cite{Ov2}.

\item dimension four: $\RR H_3$ the trivial extension of the Heisenberg Lie group, with the metric given above (\ref{e3}). 

\item dimension five: The Lie group with Lie algebra $\aa(1,2)$ and the bi-invariant metric
induced by $Q$ as in (\ref{e4}).

\item dimension six:  The Lie group with Lie algebra $\ct^* \hh_3$ and the bi-invariant metric induced by $Q$ as in (\ref{e5}).
\end{enumerate}

In the cases (3) and (4) above the corresponding Lie group  $N$ admits a cocompact lattice $\Gamma$ hence $N/\Gamma$ becomes naturally reductive.


\begin{theacknowledgments}
The author deeply thanks V. Bangert, Mathematisches Institut der Universit\"at Freiburg, for the support during the stay of the author, where part of this work was done. Many thanks to the organizers of the Meeting in Porto, for the finantial help to attend the conference.
\end{theacknowledgments}



\bibliographystyle{aipproc}   

\bibliography{sample}




\end{document}